\documentclass[12pt]{article}
\usepackage{graphics}
\usepackage{amssymb}
\usepackage{amsmath}
\usepackage{makeidx}
\usepackage{color}
\usepackage{graphicx}
\usepackage{subfigure}
\usepackage{multicol}
\usepackage{multirow}
\usepackage{float}
\usepackage{booktabs}
\usepackage[labelsep=period]{caption}
\usepackage{cite} 

\newtheorem{thm}{Theorem}[section]

\newtheorem{lem}{Lemma}[section]

\newtheorem{examp}{Example}[section]
\newtheorem{que}{Question}[section]
\numberwithin{equation}{section}

\allowdisplaybreaks

\title{\bf  Clunie lemma in several complex variables and application in PDEs
}
\author{  Wenjie Hao$^{1}$ and Qingcai Zhang$^{2}$ \\
\small $^{1}$Institute of Mathematics, Henan Academy of Sciences, Zhengzhou 450046, P.R. China; \\
\small $^{2}$School of Mathematics, Renmin University of China, Beijing 100872, P.R. China; \\
\small E-mails: haowenjie@hnas.ac.cn and zhangqcrd@163.com
}

\begin{document}
\date{}
\maketitle

\begin{quote}
{\bf Abstract:} Two purposes will be shown in this paper. The first one is to extend the classic Tumura-Clunie type theorem for meromorphic functions of one complex variable to meromorphic functions of several complex variables by using Clunie
lemma. The second one is to characterize entire solutions of certain partial differential equations in $\mathbb{C}^{m}$. Our results are extensions and generalizations of the previous theorems by Liao-Ye \cite{Liao-Ye} and Li \cite{Li11}.



\end{quote}

\section{Introduction and Main Results}
\qquad Let $f$ be a nonconstant meromorphic function in $\mathbb{C}^{m}$. We will assume that the readers are familiar with the basic notations of Nevanlinna theory of meromorphic functions in $\mathbb{C}^{m}$, such as the characteristic function $T(r,f)$ and the proximity function $m(r,f)$ (see\cite{Hayman2,Stoll5,Hu6}).

We denote by $S(r,f)$ any quantity satisfying
$$\|\quad S(r,f)=o(T(r,f)),$$
where the symbol $\|$ means that the relation holds outside a set of $r$ of finite linear measure. A meromorphic function $a$ is said to be a small function with respect to $f$ if only if $\|\quad T(r,a)=S(r,f)$. Let $\mathbb{Z}_{+}$ denote the
set of nonnegative integers. For $z=(z_{1},...,z_{m})\in\mathbb{C}^{m}$, $\mathbf{i}=(i_{1},...,i_{m})\in\mathbf{Z}_{+}^{m}$, we write
$$\partial_{z_{i}}=\frac{\partial}{\partial z_{i}},\quad (\partial_{z_{i}})^{k}=\frac{\partial^{k}}{\partial z_{i}^{k}}, \quad \partial_{z_{i}}\partial_{z_{j}}=\frac{\partial^{2}}{\partial z_{i}\partial z_{j}}, \quad i,j=1,...,m,$$
$$\partial^{\mathbf{i}}=\partial_{z}^{\mathbf{i}}=\partial_{z_{1}}^{i_{1}}\cdots\partial_{z_{m}}^{i_{m}}, \quad |\mathbf{i}|=i_{1}+\cdots+i_{m}.$$

We should pay more attention to meromorphic functions $h=P(f,\partial^{\mathbf{i}_{1}}f,...,\partial^{\mathbf{i}_{k}}f)$ which are polynomials in $f$ and the partial derivatives $\partial^{\mathbf{i}_{1}}f,..., \partial^{\mathbf{i}_{k}}f$ of $f$
with meromorphic coefficients being small functions of $f$. Such functions $h$ are called differential polynomials in $f$. The degree of the polynomial $P(x_{0},x_{1},...,x_{k})$ is called the degree  of $h$.

In 1937, Tumura\cite{Tumura9} proved the following well-known theorem.
\begin{thm}\label{thm1.1}
{\rm(see \cite{Tumura9})} Suppose that $f$ is a meromorphic function in the complex plane and has only a finite number of poles in the plane, and that $f$, $f^{(l)}$ have only a finite number of zeros for some $l\geq2$. Then
$$f=\frac{P_{1}}{P_{2}}e^{P_{3}},$$
where $P_{1}$, $P_{2}$, $P_{3}$, are polynomials. If further, $f$ and $f^{(l)}$ have no zeros, then $f(z)=e^{Az+B}$ or $f(z)=(Az+B)^{-n}$, where $A$, $B$ are constants such that $A\neq0$, and $n$ is a positive integer.
\end{thm}

However, Tumura's proof included gaps. In 1962, Clunie\cite{Clunie} presented a valid proof of the above theorem by using the following Clunie's lemma.
\begin{thm}\label{thm1.22}
{\rm(see \cite{Clunie}, Lemma 2)} Suppose that $f$ is meromorphic and transcendental in the complex plane and that
\begin{equation}\label{1.22}
f^{n}P(z,f)=Q(z,f),
\end{equation}
where $P(z,f)$, $Q(z,f)$ are differential polynomials in $f$ and the degree of $Q(z,f)$ is at most $n$. Then
$$m(r,P)=S(r,f).$$
\end{thm}

In 1964, Hayman \cite{Hayman2} proved a generalization of previous results of Tumura\cite{Tumura9} and Clunie\cite{Clunie}.
\begin{thm}\label{thm1.2}
{\rm(see \cite{Hayman2}, Theorem 3.9)} Suppose that $f$, $g$ are nonconstant meromorphic functions in the complex plane, that
\begin{equation}\label{1.2}
f^{n}+P_{n-1}(f)=g,
\end{equation}
where $P_{n-1}(f)$ denotes a differential polynomial in $f$ of degree $d\leq n-1$. If
$$N(r,f)+N\left(r,\frac{1}{g}\right)=S(r,f),$$
then $g(z)=(f(z)+\gamma)^{n}$, where $\gamma$ is meromorphic and a small function of $f$.
\end{thm}

Both Theorem {\rm\ref{thm1.1}} and Theorem {\rm\ref{thm1.2}} are shown in $\mathbb{C}^{m}$ when $m=1$. When $m\geq1$, Hu and Yang \cite{Hu3} proved the following Theorem {\rm\ref{thm1.3}}.
\begin{thm}\label{thm1.3}
{\rm(see \cite{Hu3}, Theorem 1.2)} Suppose that $f$ is meromorphic and not constant in $\mathbb{C}^{m}$, that
\begin{equation}\label{1.3}
f^{n}+P_{n-1}(f)=g,
\end{equation}
where $P_{n-1}(f)$ is a differential polynomial of degree at most $n-1$ in $f$, and that
$$\| \quad N(r,f)+N\left(r,\frac{1}{g}\right)=S(r,f).$$
Then
$$g=\left(f+\frac{a}{n}\right)^{n},$$
where $a$ is meromorphic and a small function of $f$ in $\mathbb{C}^{m}$ determined by the terms of degree $n-1$ in $P_{n-1}(f)$ and by $g$.
\end{thm}

It is natural to ask what will happen if the dominant term $f^{n}$ is replaced by $f^{n}\partial_{z_{i}}f$ in (\ref{1.3}). In fact, when $m=1$, the nonlinear differential equation
\begin{equation}\label{1.44}
f^{n}f'+P_{n-1}(f)=g,
\end{equation}
where $P_{n-1}(f)$ is a differential polynomial of degree at most $n-1$ in $f$, has attracted the attention of many scholars. ``g'', as a special function, has already been studied in early 2014 by Liao-Ye \cite{Liao-Ye}. In recent years, there has been a lot of research on equation (\ref{1.44}) (see \cite{minfeng,wei,nan}).
\begin{thm}\label{thm1.5}
{\rm(see \cite{Liao-Ye}, Theorem 1.1)} Let $P_{n-1}(f)$ be a differential polynomial in $f$ of degree $d$ with rational function coefficients. Suppose that $u$ is a nonzero rational function and $v$ is a nonconstant polynomial. If $d\leq n-1$ and
the differential equation
\begin{equation}\label{1.5}
f^{n}f'+P_{n-1}(f)=u(z)e^{v(z)}
\end{equation}
admits a meromorphic solution $f$ with finitely many poles, then $f$ has the following form:
$$f(z)=s(z)e^{v(z)/(n+1)} \quad and \quad P_{n-1}(f)\equiv0,$$
where $s(z)$ is a rational function with $s^{n}((n+1)s'+v's)=(n+1)u$. In particular, if $u$ is a polynomial, then $s$ is a polynomial, too.
\end{thm}

when $m\geq1$ is considered, the following Theorem {\rm\ref{thm1.4}} is obtained in this paper.
\begin{thm}\label{thm1.4}
Suppose that $f$ is meromorphic and not constant in $\mathbb{C}^{m}$, that
\begin{equation}\label{1.4}
f^{n}\partial_{z_{i}}f+P_{n-1}(f)=g,
\end{equation}
where $P_{n-1}(f)$ is a differential polynomial of degree at most $n-1$ in $f$, and that
$$\| \quad N(r,f)+N\left(r,\frac{1}{g}\right)=S(r,f),$$
$$\| \quad T(r,f(z_{[i]}))\neq S(r,f),$$
where $z_{[i]}=(\zeta_{1},...,\zeta_{i-1},z_{i},\zeta_{i+1},...,\zeta_{m})$, $\zeta_{1},...,\zeta_{m}$ are any fixed complex numbers.
Then
$$g=af^{n+1}, \quad P_{n-1}(f)\equiv0,$$
where $a$ is meromorphic and a small function of $f$ in $\mathbb{C}^{m}$ determined by the terms of degree $n-1$ in $P_{n-1}(f)$ and by $g$.
\end{thm}

Obviously, Theorem {\rm\ref{thm1.5}} is a special case of Theorem {\rm\ref{thm1.4}}. But it must be indicated that the proofs of these two theorems are slightly different. For example, the author uses the multiplicity of zeros on both sides of
the equation in Theorem {\rm\ref{thm1.5}}. This is not true in several complex variables, since the multiplicity of a zero of $f(z_{1}, z_{2},..., z_{n})-a$ might be lower than or equal to the one for its partial derivative $\partial_{z_{i}}f$.
Another difference lies in the integration of variables. The constant calculated by integration in one complex variable becomes the function which is worked out  in several complex variables.  Therefore, it can be believed that Theorem
{\rm\ref{thm1.5}} has been improved in a different way. Now, we give two examples to illustrate the necessity of the conditions in Theorem {\rm\ref{thm1.4}}.

\begin{examp}\label{examp1.1}
The function $f(z_{1},z_{2})=e^{z_{1}}+z_{2}$ is a transcendental solution of the differential equation
$$f^{2}\partial_{z_{2}}f-2z_{2}\partial_{z_{1}}f-z_{2}^{2}=e^{2z_{1}}.$$
\end{examp}

{\rm
By Example \ref{examp1.1}, we can observe that for $T(r,f(z_{[i]}))=T(r,f(\zeta_{1},z_{2}))=S(r,f)$, where $\zeta_{1}$ is a complex number, $P_{n-1}(f)=-2z_{2}\partial_{z_{1}}f-z_{2}^{2}\not\equiv0$. Therefore, the condition
``$T(r,f(z_{[i]}))\neq S(r,f)$" in Theorem \ref{thm1.4} is necessary.

\begin{examp}\label{examp1.2}
The function $f(z_{1},z_{2})=e^{z_{1}+z_{2}}+z_{1}+z_{2}$ is a transcendental solution of the differential equation
$$f^{2}\partial_{z_{i}}f-[2(z_{1}+z_{2})+1](\partial_{z_{i}}f)^{2}-[(z_{1}+z_{2})^{2}-2(z_{1}+z_{2})-2]\partial_{z_{i}}f-1=e^{3(z_{1}+z_{2})}.$$
\end{examp}

{\rm
By Example \ref{examp1.2}, we can observe that for $n=d=2$, $P_{n-1}(f)=-[2(z_{1}+z_{2})+1]\partial_{z_{i}}f^{2}-[(z_{1}+z_{2})^{2}-2(z_{1}+z_{2})-2]\partial_{z_{i}}f-1\not\equiv0$. Therefore, the condition ``$d\leq n-1$" in Theorem \ref{thm1.4}
is necessary.

In view of Theorem \ref{thm1.4}, we turn our attention to entire solutions of general inviscid Burgers' equation. Burgers' equation or Bateman-Burgers equation
\begin{equation}\label{1.18}
f\partial_{z_{2}}f+\partial_{z_{1}}f=\beta(\partial_{z_{2}})^{2}f, \quad f=f(z_{1},z_{2}),
\end{equation}
where $\beta\neq0$ is a constant, is a basic partial differential equation in modeling many physical phenomena, such as fluid mechanics, traffic flow, growth of interfaces, and financial mathematics. This is the simplest partial differential
equation which combines both nonlinear propagation effects and diffusive effects. When the right term of (\ref{1.18}) is removed, one gets the hyperbolic partial differential equation
$$f\partial_{z_{2}}f+\partial_{z_{1}}f=0, \quad f=f(z_{1},z_{2}),$$
which is called inviscid Burgers' equation. In 2005, Li considered entire solutions of inviscid Burgers' equation and obtained the following theorem.
\begin{thm}\label{thm1.19}
{\rm(see \cite{Li11}, Theorem 2.3)} A function $f$ is an entire solution of the partial differential equation
\begin{equation}\label{1.19}
cf^{n}\partial_{z_{2}}f+\partial_{z_{1}}f=0
\end{equation}
in $\mathbb{C}^{2}$, where $c\neq0$ is a constant and $n\geq0$ is an integer, if and only if $f$ is a constant when $n>0$;
and $f=g(x+ct)$ when $n=0$, where $g$ is an entire function in the complex plane.
\end{thm}

Recently, L${\rm\ddot{u}}$ \cite{LF13} has changed the right term of (\ref{1.19}) to a polynomial and deduced the entire solution of general inviscid Burgers' equation combining characteristic equations for quasi-linear partial differential
equations with normal family and the Nevanlinna theory. More generally, we obtain the following theorem.
\begin{thm}\label{thm1.20}
Suppose that $f$ is entire and not constant in $\mathbb{C}^{m}$, that
\begin{equation}\label{1.20}
f^{n}\partial_{z_{i}}f+P_{n-1}(f)=g(f),
\end{equation}
where $P_{n-1}(f)$ is a differential polynomial of degree $d$ at most $n-1$ in $f$, $g$ is a meromorphic function in $\mathbb{C}$.
Then $f$ must be one of the following form:
\vskip0.1in
{\rm(1)} $f=az_{i}+\alpha(z_{1},...,z_{i-1},z_{i+1},...,z_{m})$,
\vskip0.1in
{\rm(2)} $f=a+e^{bz_{i}+\alpha(z_{1},...,z_{i-1},z_{i+1},...,z_{m})}$,
\vskip0.1in
where $a$, $b$ are complex numbers, and $\alpha(z_{1},...,z_{i-1},z_{i+1},...,z_{m})$ is an entire function in $\mathbb{C}^{m-1}$.
\end{thm}

We also give examples to show that each form of entire solutions in Theorem \ref{thm1.20} can indeed occur. Let $\alpha(z_{2},z_{3},...,z_{m-1})$ be an arbitrary entire function in $\mathbb{C}^{m-1}$ in all the following examples.

\begin{examp}\label{examp1.3}
The function $f(z_{1},z_{2})=z_{1}+z_{2}$ is an entire solution of the partial differential equation
$$f^{2}\partial_{z_{i}}f+\partial_{z_{i}}f=f^{2}+1,$$
and is the form {\rm(1)} in Theorem {\rm\ref{thm1.20}}.
\end{examp}

\begin{examp}\label{examp1.4}
The function $f(z_{1},z_{2})=e^{z_{1}+z_{2}}$ is an entire solution of the partial differential equation
$$f^{2}\partial_{z_{i}}f-f=f^{2}-f,$$
and is the form {\rm(2)} in Theorem {\rm\ref{thm1.20}}.
\end{examp}

At the end of this section, we pose a question.
\begin{que}\label{que1}
The conclusion is still true for ``$d\geq n$" in Theorem {\rm\ref{thm1.20}} {\rm?}
\end{que}

Here's an example to answer this question, but we have no idea to solve it.

\begin{examp}\label{examp1.5}
The function $f(z_{1},z_{2})=e^{z_{1}+z_{2}}$ is an entire solution of the partial differential equation
$$f^{n}\partial_{z_{i}}f-(\partial_{z_{i}}f)^{k}=f^{n+1}-f^{k},$$
where $k\geq n$, the solution $f$ is of the form {\rm(2)} in Theorem {\rm\ref{thm1.20}}.
\end{examp}

\section{Some lemmas}

\qquad To prove our results, we will employ Nevanlinna theory. For the convenience of the reader who is not familiar with Nevanlinna theory, we list here the notation and results from Nevanlinna theory which will be used in the proofs (see e.g.,
\cite{Stoll5,Hu6,Vitter14,Griffiths16,Noguchi17,Ru18}).

Set $\|z\|=(|z_{1}|^{2}+\cdots+|z_{m}|^{2})^{1/2}$ for $z=(z_{1},...,z_{m})\in\mathbb{C}^{m}$. For $r>0$, define
$$B_{m}(r)=\{z\in\mathbb{C}^{m}:\|z\|<r\},\quad S_{m}(r)=\{z\in\mathbb{C}^{m}:\|z\|=r\}.$$
Let $d=\partial+\overline{\partial}$, $d^{c}=(4\pi\sqrt{-1})^{-1}(\partial-\overline{\partial})$. Thus, $dd^{c}=(2\pi)^{-1}\sqrt{-1}\partial\overline{\partial}$. Write
$$\upsilon_{m}(z)=(dd^{c}\|z\|^{2})^{m-1},\quad \sigma_{m}(z)=d^{c}\log\|z\|^{2}\wedge(dd^{c}\log\|z\|^{2})^{m-1},$$
for $z\in\mathbb{C}^{m}\backslash\{0\}$.

For a divisor $\nu$ on $\mathbb{C}^{m}$, define the following counting function of $\nu$ by
$$N(r,\nu)=\int_{1}^{r}\frac{n(t)}{t^{2m-1}}dt,\quad(1<r<\infty),$$
where
$$n(t)=\left\{
\begin{array}{rcl}
\int_{|\nu|\cap B(t)}\nu(z)\upsilon_{m}(z),& &{m\geq2},\\
\sum_{|z|\leq t}\nu(z),& &{m=1.}
\end{array}\right.$$

Let $f(z)$ be a nonzero entire function on $\mathbb{C}^{m}$. For a point $z_{0}\in\mathbb{C}^{m}$, we write$f(z)=\sum_{i=0}^{\infty}P_{i}(z-z_{0})$, where the term $P_{i}$ is a homogeneous polynomial of degree $i$. Denote the zero-multiplicity of
$f$ at $z_{0}$ by $\nu_{f}(z_{0})=\min\{i:P_{i}\neq0\}$ in the sense of \cite{Fujimoto15}. Set $|\nu_{f}|={\rm Supp}\nu_{f}$, which is a pure $(n-1)$-dimensional analytic subset or empty set.

Let $f(z)$ be a nonzero meromorphic function on $\mathbb{C}^{m}$ in the sense that f can be written as a quotient of two relatively prime holomorphic functions. For each $z\in\mathbb{C}^{m}$, write $f(z)=(f_{1}(z), f_{2}(z))$ where
$f_{1}(\not\equiv0)$, $f_{2}$ are two relatively prime holomorphic functions such that dim$\{z\in\mathbb{C}^{m}:f_{1}(z)=f_{2}(z)=0\}\leq m-2$. Thus, $f$ may be regarded as a meromorphic mapping from $\mathbb{C}^{m}$ into $\mathbb{P}^{1}$ such
that $f^{-1}(\infty)\neq\mathbb{C}^{m}$.

Define $\nu_{f}^{0}=\nu_{f_{2}}$, $\nu_{f}^{\infty}=\nu_{f_{1}}$, which are independent of the choice of $f_{1}$, $f_{2}$. For a meromorphic function $f$ on $\mathbb{C}^{m}$, we usually write $N(r,f)=N(r,\nu_{f}^{\infty})$ and
$N\left(r,\frac{1}{f}\right)=N(r,\nu_{f}^{0})$. Denote by $m(r,f)$ the proximity function of $f$ defined as
$$m(r,f)=\int_{S_{m}(r)}\log^{+}\max\{|f(z)|\}\sigma_{m}(z),$$
where $\log^{+}x=\max\{\log x,0\}$. Then the Nevanlinna characteristic function of $f$ is defined as $T(r,f)=N(r,f)+m(r,f)$. Then the first main theorem is said that
$$T\left(r,\frac{1}{f-a}\right)=T(r,f)+O(1),$$
for any value $a\in\mathbb{c}\cup\{\infty\}$. The Cartan characteristic function is defined by
$$T_{f}(r)=\int_{S_{m}(r)}\log\max\{|f_{1}(z)|,|f_{2}(z)|\}\sigma_{m}(z)-\int_{S_{m}(1)}\log\max\{|f_{1}(z)|,|f_{2}(z)|\}\sigma_{m}(z).$$
The two characteristic functions have the relation $T_{f}(r)=T(r,f)+O(1)$.
The order $\rho(f)$ of a meromorphic function $f$ is defined as follows:
$$\rho\left(f\right)=\limsup_{r\to\infty}\frac{\log T\left(r,f\right)}{\log r}.$$

In order to prove the main theorems in this paper, the following auxiliary lemmas from Nevanlinna theory in $\mathbb{C}^{m}$ are needed. The following form of the lemma of logarithmic derivative has been frequently used in value distribution
theory of meromorphic functions in $\mathbb{C}^{m}$.

\begin{lem}\label{lem2.11}
{\rm(see \cite{Hu6,Vitter14})} Assume that $f$ is a nonconstant meromorphic function in $\mathbb{C}^{m}$. Let $\nu=(\nu_{1},...,\nu_{m})\in\mathbb{Z}^{m}_{+}$ be a multi-index. Then for any $\varepsilon>0$,
$$\|\quad m\left(r,\frac{\partial^{\nu}f}{f}\right)\leq|\nu|\log^{+}T(r,f)+|\nu|(1+\varepsilon)\log^{+}\log T(r,f)+O(1).$$
\end{lem}

Proofs of Theorem {\rm\ref{thm1.4}} and Theorem {\rm\ref{thm1.20}} are based on the following generalised Clunie lemma in $\mathbb{C}^{m}$. It has a more general form than Theorem {\rm\ref{thm1.22}}. For the case $m=1$, see \cite{He-Xiao12}. For
some special cases, see \cite{Hu6,Hu7}. A general proof can be found in \cite{Hu8}.
\begin{lem}\label{lem2.1}
{\rm(see \cite{Hu3}, Lemma 2.1)} Let $f$ be a nonconstant meromorphic function on $\mathbb{C}^{m}$. Take a positive integer $n$ and take polynomials of $f$ and its partial derivatives:
\begin{equation*}
P(f)=\sum_{\mathbf{p}\in I}a_{\mathbf{p}}f^{p_{0}}(\partial^{\mathbf{i}_{1}}f)^{p_{1}}\cdots (\partial^{\mathbf{i}_{l}}f)^{p_{l}}, \mathbf{p}=(p_{0},...,p_{l})\in\mathbb{Z}_{+}^{l+1},
\end{equation*}
\begin{equation*}
Q(f)=\sum_{\mathbf{q}\in J}c_{\mathbf{q}}f^{q_{0}}(\partial^{\mathbf{j}_{1}}f)^{q_{1}}\cdots (\partial^{\mathbf{j}_{s}}f)^{q_{s}}, \mathbf{q}=(q_{0},...,q_{l})\in\mathbb{Z}_{+}^{s+1},
\end{equation*}
and
\begin{equation*}
B(f)=\sum_{k=0}^{n}b_{k}f^{k},
\end{equation*}
where $I$, $J$ are finite sets of distinct elements and $a_{\mathbf{p}}$, $c_{\mathbf{q}}$, $b_{k}$ are meromorphic functions on $\mathbb{C}^{m}$ with $b_{n}\not\equiv0$. Assume that f satisfies the equation
$$B(f)Q(f)=P(f)$$
such that $P(f)$, $Q(f)$ and $B(f)$ are differential polynomials, that is, their coefficients being small functions of $f$. If $\deg(P(f))\leq n= \deg(B(f))$, then
$$\| \quad m(r,Q(f))=S(r,f).$$
\end{lem}

\begin{lem}\label{lem2.2}
{\rm(see \cite{Li10}, Theorem 1)} Let $F$ be an entire solution in $\mathbb{C}^{m}$ of the partial differential equation
$\partial w/\partial z_{j}=f(w)$, say $j=1$, where $f$ is a meromorphic function in $\mathbb{C}$. Then $F$ must be of
the form $F=az_{1}+\alpha(z_{2},z_{3},...,z_{m})$ or of the form $F=a+e^{bz_{1}+\alpha(z_{2},z_{3},...,z_{m})}$, where $a$, $b$ are complex numbers, and $\alpha(z_{2},z_{3},...,z_{m})$ is an entire function in $\mathbb{C}^{m-1}$.
\end{lem}

\begin{lem}\label{lem2.4}
{\rm(see \cite{Hu6}, Theorem 1.26)} Let $f(z)$ be a nonconstant meromorphic function in the parabolic manifold $M$. Assume that
$$R(z,w)=\frac{A(z,w)}{B(z,w)}.$$
Then
$$T(r,R_{f})=\max\{p, q\}T(r,f)+O\left(\sum_{j=0}^{p}T(r,a_{j})+\sum_{j=0}^{q}T(r,b_{j})\right),$$
where $R_{f}(z)=R(z,f(z))$ and two relatively prime polynomials $A(z,w)$, $B(z,w)$ are given respectively as follows:
$$A(z,w)=\sum_{j=0}^{p}a_{j}w^{j},\quad B(z,w)=\sum_{j=0}^{q}b_{j}w^{j}.$$
\end{lem}

\begin{lem}\label{lem2.5}
{\rm(see \cite{Chang-Li-Yang19}, Theorem 4.1)} If $g$ is a transcendental meromorphic function in $\mathbb{C}$ and $f$ is a transcendental entire function in $\mathbb{C}^{m}$, then
$$\lim_{r\rightarrow\infty}\frac{T(r,g(f))}{T(r,f)}=\infty.$$
\end{lem}

\begin{lem}\label{lem2.6}
{\rm(see \cite{Vitter14})} If $f$ is a meromorphic function in $\mathbb{C}^{m}$, then
$$\|\quad T(r,\partial_{z_{i}}f)=O(T(r,f)).$$
\end{lem}

\section{Proof of theorem 1.6}
\quad Set $P=P_{n-1}(f)$. Without loss of generality, we say $i=1$ in (\ref{1.4}). We now differentiate (\ref{1.4}) and obtain
\begin{equation}\label{3.1}
f^{n}\partial_{z_{1}}\partial_{z_{k}}f+nf^{n-1}\partial_{z_{1}}f\partial_{z_{k}}f+\partial_{z_{k}}P=\partial_{z_{k}}g,\quad  k=1,2,...,m,
\end{equation}
where $\partial_{z_{k}}P$ is a differential polynomial in $f$ of degree at most $n-1$. We multiply (\ref{1.4}) by $\frac{\partial_{z_{k}}g}{g}$ and subtract from (\ref{3.1}). This gives
\begin{equation}\label{3.2}
f^{n-1}\varphi_{k}=\partial_{z_{k}}P-\frac{\partial_{z_{k}}g}{g}P,
\end{equation}
where
\begin{equation}\label{3.3}
\varphi_{k}=\frac{\partial_{z_{k}}g}{g}f\partial_{z_{1}}f-f\partial_{z_{1}}\partial_{z_{k}}f-n\partial_{z_{1}}f\partial_{z_{k}}f.
\end{equation}
Note that by (\ref{1.4}) $g$ has poles only at poles of $f$ or of the coefficients of $P$ and so we have
\begin{equation*}
\|\quad N(r,g)=S(r,f).
\end{equation*}
Thus
\begin{equation*}
\|\quad N\left(r,\frac{\partial_{z_{k}}g}{g}\right)\leq N(r,g)+N\left(r,\frac{1}{g}\right)=S(r,f)
\end{equation*}
by hypothesis. Apply Lemmas \ref{lem2.4} and \ref{lem2.6} to (\ref{1.4}), so that
$$\|\quad T(r,g)=O(T(r,f)).$$
By using the lemma of logarithmic derivative, we know
\begin{equation*}
\|\quad m\left(r,\frac{\partial_{z_{k}}g}{g}\right)=S(r,g)=S(r,f)
\end{equation*}
and so
\begin{equation*}
\|\quad T\left(r,\frac{\partial_{z_{k}}g}{g}\right)=m\left(r,\frac{\partial_{z_{k}}g}{g}\right)+N\left(r,\frac{\partial_{z_{k}}g}{g}\right)=S(r,f).
\end{equation*}
Thus (\ref{3.2}) satisfies the conditions of the generalised Clunie lemma, we deduce that
\begin{equation*}
\|\quad m(r,\varphi_{k})=S(r,f).
\end{equation*}
Note that $\varphi_{k}$ has poles only at poles of $f$ or $\frac{\partial_{z_{k}}g}{g}$ and so we have
\begin{equation*}
\|\quad N(r,\varphi_{k})\leq4N(r,f)+N\left(r,\frac{\partial_{z_{k}}g}{g}\right)=S(r,f),
\end{equation*}
and hence
\begin{equation*}
\|\quad T(r,\varphi_{k})=m(r,\varphi_{k})+N(r,\varphi_{k})=S(r,f),
\end{equation*}
for each $k=1,...,m$. Suppose that $P\not\equiv0$. We distinguish three cases to discuss.

Case 1: $\varphi_{1}\equiv0$ and $\varphi_{k}\equiv0$ for any $k=2,3,...,m$ holds. By (\ref{3.2}), we have
\begin{equation*}
\frac{\partial_{z_{k}}g}{g}=\frac{\partial_{z_{k}}P}{P}, \quad g=cP, \quad k=1,2,3,...,m,
\end{equation*}
where $c$ is a nonzero constant. Substituting it into (\ref{1.4}) yields
\begin{equation}\label{3.4.1}
f^{n}\partial_{z_{1}}f=(c-1)P.
\end{equation}

If $c=1$, we have $f\equiv0$ or $\partial_{z_{1}}f\equiv0$. This contradicts the condition $T(r,f(z'))\neq S(r,f)$.

If $c\neq1$, (\ref{3.4.1}) satisfies the conditions of the generalised Clunie lemma, we obtain
$$\|\quad m(r,\partial_{z_{1}}f)=S(r,f), \quad m(r,f\partial_{z_{1}}f)=S(r,f).$$
Moreover,
$$\|\quad T(r,\partial_{z_{1}}f)=S(r,f), \quad T(r,f\partial_{z_{1}}f)=S(r,f),$$
by hypothesis. Thus,
$$\|\quad T(r,f)\leq T(r,f\partial_{z_{1}}f)+T\left(r,\frac{1}{\partial_{z_{1}}f}\right)=S(r,f),$$
which yields a contradiction.

Case 2: $\varphi_{1}\equiv0$ and $\varphi_{k}\not\equiv0$ for some $k=2,3,...,m$ holds. By (\ref{3.3}), we have
\begin{equation*}\label{3.4}
\frac{\partial_{z_{1}}g}{g}=\frac{(\partial_{z_{1}})^{2}f}{\partial_{z_{1}}f}+n\frac{\partial_{z_{1}}f}{f}.
\end{equation*}
Integrating this equality yields that
$$f^{n}\partial_{z_{1}}f=ge^{c(\widehat{z_{1}})},$$
where $c(\widehat{z_{1}})$ is a entire function in $(z_{2},z_{3},...,z_{m})$. Hence,
\begin{equation}\label{3.19}
N\left(r,\frac{1}{f}\right)\leq N\left(r,\frac{1}{g}\right)=S(r,f).
\end{equation}
On the other hand, by (\ref{3.3}), we get
$$\frac{1}{f^{2}}=\frac{1}{\varphi_{k}}\left(\frac{\partial_{z_{k}}g}{g}\frac{\partial_{z_{1}}f}{f}-\frac{\partial_{z_{1}}\partial_{z_{k}}f}{f}
-n\frac{\partial_{z_{1}}f}{f}\frac{\partial_{z_{k}}f}{f}\right).$$
By the logarithmic derivative lemma, we have
\begin{equation}\label{3.20}
2m\left(r,\frac{1}{f}\right)\leq m\left(r,\frac{1}{\varphi_{k}}\right)+S(r,f)\leq T(r,\varphi_{k})+S(r,f)=S(r,f),
\end{equation}
combining (\ref{3.19}) with (\ref{3.20}) yields that
$$T(r,f)=T\left(r,\frac{1}{f}\right)+O(1)=m\left(r,\frac{1}{f}\right)+N\left(r,\frac{1}{f}\right)+O(1)=S(r,f).$$
This is absurd.

Case 3: $\varphi_{1}\not\equiv0$. We differentiate $\varphi_{1}$ and obtain
\begin{align}\label{3.5}
\partial_{z_{1}}\varphi_{1}&=f\left[\frac{g(\partial_{z_{1}})^{2}g
-(\partial_{z_{1}}g)^{2}}{g^{2}}\partial_{z_{1}}f+\frac{\partial_{z_{1}}g}{g}(\partial_{z_{1}})^{2}f-(\partial_{z_{1}})^{3}f\right]\nonumber\\
&+\partial_{z_{1}}f\left[\frac{\partial_{z_{1}}g}{g}\partial_{z_{1}}f-(2n+1)(\partial_{z_{1}})^{2}f\right].
\end{align}
From (\ref{3.3}) and (\ref{3.5}), we have
\begin{equation}\label{3.6}
f[B_{1}\partial_{z_{1}}f+B_{2}(\partial_{z_{1}})^{2}f-(\partial_{z_{1}})^{3}f]+\partial_{z_{1}}f[A\partial_{z_{1}}f-(2n+1)(\partial_{z_{1}})^{2}f]=0,
\end{equation}
where
\begin{align*}
B_{1}&=\frac{g(\partial_{z_{1}})^{2}g-(\partial_{z_{1}}g)^{2}}{g^{2}}-\frac{\partial_{z_{1}}\varphi_{1}}{\varphi_{1}}\frac{\partial_{z_{1}}g}{g},\nonumber\\
B_{2}&=\frac{\partial_{z_{1}}\varphi_{1}}{\varphi_{1}}+\frac{\partial_{z_{1}}g}{g},\nonumber\\
A&=\frac{n\partial_{z_{1}}\varphi_{1}}{\varphi_{1}}+\frac{\partial_{z_{1}}g}{g}.
\end{align*}
If $z_{0}$ is the common zero of $f$ and $\partial_{z_{1}}f$, but is not a zero or pole of $\frac{\partial_{z_{1}}g}{g}$, then from (\ref{3.3}) we have $z_{0}$ is a zero of $\varphi_{1}$. Suppose that $z_{0}$ is a zero of $f$, but is not a zero
or pole of $\partial_{z_{1}}f$ and $\frac{\partial_{z_{1}}g}{g}$. Then (\ref{3.6}) gives $A\partial_{z_{1}}f(z_{0})-(2n+1)(\partial_{z_{1}})^{2}f(z_{0})=0$. Set
\begin{equation}\label{3.60}
Q=\frac{A\partial_{z_{1}}f-(2n+1)(\partial_{z_{1}})^{2}f}{f}
=\frac{B_{1}\partial_{z_{1}}f+B_{2}(\partial_{z_{1}})^{2}f-(\partial_{z_{1}})^{3}f}{\partial_{z_{1}}f}.
\end{equation}
Note that $N(r,Q)=S(r,f)$, and by using the lemma of logarithmic derivative, we know
$$\|\quad T(r,Q)=m(r,Q)+N(r,Q)=S(r,f).$$

We distinguish two subcases to discuss.

Subcase 3.1: $Q\equiv0$. Then
\begin{equation}\label{3.40}
\frac{(\partial_{z_{1}})^{2}f}{\partial_{z_{1}}f}=\frac{A}{(2n+1)}
=\frac{n\partial_{z_{1}}\varphi_{1}}{(2n+1)\varphi_{1}}+\frac{\partial_{z_{1}}g}{(2n+1)g},
\end{equation}
that gives
\begin{equation}\label{3.8}
\partial_{z_{1}}f=\varphi_{1}^{\frac{n}{2n+1}}g^{\frac{1}{2n+1}}e^{c(\widehat{z_{1}})},
\end{equation}
upon integrating (\ref{3.40}), where $c(\widehat{z_{1}})$ is a entire function in $(z_{2},z_{3},...,z_{m})$. Note that $(\partial_{z_{1}})^{2}f=\frac{A}{(2n+1)}\partial_{z_{1}}f$. Substituting (\ref{3.8}) into (\ref{3.3}), we get
\begin{align}\label{3.9.1}
\varphi_{1}&=\left(\frac{\partial_{z_{1}}g}{g}-\frac{A}{2n+1}\right)f\partial_{z_{1}}f-n(\partial_{z_{1}}f)^{2}\nonumber\\
&=\frac{n}{2n+1}\left(\frac{2\partial_{z_{1}}g}{g}-\frac{\partial_{z_{1}}\varphi_{1}}{\varphi_{1}}\right)f\varphi_{1}^{\frac{n}{2n+1}}g^{\frac{1}{2n+1}}e^{c(\widehat{z_{1}})}
-n\varphi_{1}^{\frac{2n}{2n+1}}g^{\frac{2}{2n+1}}e^{2c(\widehat{z_{1}})}.
\end{align}

Now we prove
\begin{equation}\label{3.99}
\frac{2\partial_{z_{1}}g}{g}-\frac{\partial_{z_{1}}\varphi_{1}}{\varphi_{1}}\not\equiv0.
\end{equation}
Otherwise, we have
\begin{equation}\label{3.999}
g=\varphi_{1}^{\frac{1}{2}}e^{c(\widehat{z_{1}})}.
\end{equation}
We may fix complex numbers $\zeta_{2},...,\zeta_{m}$ such that $g(z_{1},\zeta_{2},...,\zeta_{m})$ is a function in $z_{1}$. By (\ref{3.999}) we have
\begin{equation*}
T(r,g(z_{1},\zeta_{2},...,\zeta_{m}))=O(T(r,\varphi_{1}(z_{1},\zeta_{2},...,\zeta_{m})))=S(r,f).
\end{equation*}
Note that (\ref{1.4}) holds for $z_{[i]}:=(z_{1},\zeta_{2},...,\zeta_{m})$, we get
\begin{equation}\label{3.991}
f(z_{[i]})^{n}\partial_{z_{1}}f(z_{[i]})+P_{n-1}(f(z_{[i]}))=g(z_{[i]}).
\end{equation}
From (\ref{1.4}), Lemmas \ref{lem2.4} and \ref{lem2.6}, we obtain
$$T(r,g(z_{[i]}))=O(T(r,f(z_{[i]}))).$$
Thus, we have $\rho(g(z_{[i]}))\leq\rho(f(z_{[i]}))$.

If $\rho(g(z_{[i]}))=\rho(f(z_{[i]}))$, then
$$T(r,f(z_{[i]}))=S(r,f),$$
which yields a contradiction.

If $\rho(g(z_{[i]}))<\rho(f(z_{[i]}))$, then
$$T(r,g(z_{[i]}))=S(r,f(z_{[i]})).$$
By using the generalised Clunie lemma, we obtain
$$\|\quad m(r,\partial_{z_{1}}f(z_{[i]}))=S(r,f(z_{[i]})), \quad m(r,f(z_{[i]})\partial_{z_{1}}f(z_{[i]}))=S(r,f(z_{[i]})),$$
and so
$$\|\quad T(r,\partial_{z_{1}}f(z_{[i]}))=S(r,f(z_{[i]})), \quad T(r,f(z_{[i]})\partial_{z_{1}}f(z_{[i]}))=S(r,f(z_{[i]}))$$
by hypothesis. Thus,
$$\|\quad T(r,f(z_{[i]}))\leq T(r,f(z_{[i]})\partial_{z_{1}}f(z_{[i]}))+T\left(r,\frac{1}{\partial_{z_{1}}f(z_{[i]})}\right)=S(r,f(z_{[i]})),$$
which also yields a contradiction.

Next we denote that
$$Y=\frac{2n+1}{\frac{2\partial_{z_{1}}g}{g}-\frac{\partial_{z_{1}}\varphi_{1}}{\varphi_{1}}}.$$
From (\ref{3.9.1}) we have
\begin{equation}\label{3.10}
f=\frac{Y}{n}\varphi_{1}^{\frac{n+1}{2n+1}}g^{-\frac{1}{2n+1}}e^{-c(\widehat{z_{1}})}
+Y\varphi_{1}^{\frac{n}{2n+1}}g^{\frac{1}{2n+1}}e^{c(\widehat{z_{1}})}.
\end{equation}
Differentiate (\ref{3.10}) and obtain
\begin{align}\label{3.11}
\partial_{z_{1}}f&=\left[\frac{\partial_{z_{1}}Y}{n}+\frac{(n+1)Y}{(2n+1)n}\frac{\partial_{z_{1}}\varphi_{1}}{\varphi_{1}}
-\frac{Y}{(2n+1)n}\frac{\partial_{z_{1}}g}{g}\right]\varphi_{1}^{\frac{n+1}{2n+1}}g^{-\frac{1}{2n+1}}e^{-c(\widehat{z_{1}})}\nonumber\\
&~+\left(\partial_{z_{1}}Y+\frac{nY}{2n+1}\frac{\partial_{z_{1}}\varphi_{1}}{\varphi_{1}}+\frac{Y}{2n+1}\frac{\partial_{z_{1}}g}{g}\right)\varphi_{1}^{\frac{n}{2n+1}}g^{\frac{1}{2n+1}}e^{c(\widehat{z_{1}})}.
\end{align}
From (\ref{3.8}) and (\ref{3.11}), we have
\begin{align}\label{3.111}
&\left(1-\frac{nY}{2n+1}\frac{\partial_{z_{1}}\varphi_{1}}{\varphi_{1}}
-\partial_{z_{1}}Y-\frac{Y}{2n+1}\frac{\partial_{z_{1}}g}{g}\right)\varphi_{1}^{\frac{-1}{2n+1}}g^{\frac{2}{2n+1}}e^{2c(\widehat{z_{1}})}\nonumber\\
&=\frac{\partial_{z_{1}}Y}{n}-\frac{Y}{(2n+1)n}\frac{\partial_{z_{1}}g}{g}+\frac{(n+1)Y}{(2n+1)n}\frac{\partial_{z_{1}}\varphi_{1}}{\varphi_{1}}.
\end{align}

If $1-\frac{nY}{2n+1}\frac{\partial_{z_{1}}\varphi_{1}}{\varphi_{1}}
-\partial_{z_{1}}Y-\frac{Y}{2n+1}\frac{\partial_{z_{1}}g}{g}\equiv0$, then $\frac{\partial_{z_{1}}Y}{n}-\frac{Y}{(2n+1)n}\frac{\partial_{z_{1}}g}{g}+\frac{(n+1)Y}{(2n+1)n}\frac{\partial_{z_{1}}\varphi_{1}}{\varphi_{1}}\equiv0$. Thus,
$$\frac{1}{2n+1}\frac{\partial_{z_{1}}g}{g}-\frac{n+1}{2n+1}\frac{\partial_{z_{1}}\varphi_{1}}{\varphi_{1}}=\frac{\partial_{z_{1}}Y}{Y}, \quad g=Y^{2n+1}\varphi_{1}^{-n-1}e^{c(\widehat{z_{1}})},$$
where $c(\widehat{z_{1}})$ is a entire function in $(z_{2},z_{3},...,z_{m})$. For the same reason as proof of (\ref{3.999}), we can get the contradiction.

If $1-\frac{nY}{2n+1}\frac{\partial_{z_{1}}\varphi_{1}}{\varphi_{1}}
-\partial_{z_{1}}Y-\frac{Y}{2n+1}\frac{\partial_{z_{1}}g}{g}\not\equiv0$, from (\ref{3.111}) we get
$$g=\omega e^{c(\widehat{z_{1}})},$$
where $\omega$ is a small function of $f$ and $c(\widehat{z_{1}})$ is a entire function in $(z_{2},z_{3},...,z_{m})$. We also can get the contradiction.

Subcase 3.2: $Q\not\equiv0$. From (\ref{3.60}) we have
\begin{equation}\label{3.7}
A\partial_{z_{1}}f-(2n+1)(\partial_{z_{1}})^{2}f-Qf=0.
\end{equation}
\begin{equation}\label{3.50}
(B_{1}-Q)\partial_{z_{1}}f+B_{2}(\partial_{z_{1}})^{2}f-(\partial_{z_{1}})^{3}f=0.
\end{equation}
We differentiate (\ref{3.7}) and obtain
\begin{equation}\label{3.12}
(\partial_{z_{1}}A-Q)\partial_{z_{1}}f+A(\partial_{z_{1}})^{2}f-(2n+1)(\partial_{z_{1}})^{3}f-f\partial_{z_{1}}Q=0.
\end{equation}
By eliminating $(\partial_{z_{1}})^{3}f$ from (\ref{3.50}) and (\ref{3.12}) we have
\begin{equation}\label{3.51}
[(2n+1)B_{1}-2nQ-\partial_{z_{1}}A]\partial_{z_{1}}f+[(2n+1)B_{2}-A](\partial_{z_{1}})^{2}f+f\partial_{z_{1}}Q=0.
\end{equation}
It follows eliminating $(\partial_{z_{1}})^{2}f$ from (\ref{3.7}) and (\ref{3.51}) that
\begin{equation*}\label{3.13}
R_{1}\partial_{z_{1}}f=R_{2}f,
\end{equation*}
where
\begin{align*}
R_{1}&=(2n+1)\partial_{z_{1}}A-(2n+1)^{2}B_{1}+2n(2n+1)Q+A^{2}-(2n+1)AB_{2},\nonumber\\
R_{2}&=Q[A-(2n+1)B_{2}]+(2n+1)\partial_{z_{1}}Q.
\end{align*}
If $R_{1}\equiv0$, then $R_{2}\equiv0$. Thus,
$$\frac{\partial_{z_{1}}Q}{Q}=B_{2}-\frac{A}{2n+1}.$$
That is,
$$\frac{\partial_{z_{1}}Q}{Q}=\frac{n+1}{2n+1}\frac{\partial_{z_{1}}\varphi_{1}}{\varphi_{1}}+\frac{2n}{2n+1}\frac{\partial_{z_{1}}g}{g}.$$
By integrating we get
$$Q=\varphi_{1}^{-\frac{n+1}{2n+1}}g^{-\frac{2n}{2n+1}}e^{c(\widehat{z_{1}})},$$
where $c(\widehat{z_{1}})$ is a entire function in $(z_{2},z_{3},...,z_{m})$. Moreover,
$$g=Q^{\frac{2n+1}{2n}}\varphi_{1}^{\frac{n+1}{2n}}e^{-\frac{2n+1}{2n}c(\widehat{z_{1}})},$$
we can get a contradiction by the same argument in (\ref{3.999}).

If $R_{1}\not\equiv0$, then
$$\partial_{z_{1}}f=V_{1}f,\quad (\partial_{z_{1}})^{2}f=V_{2}f,$$
where $V_{1}=\frac{R_{2}}{R_{1}}$ and $V_{2}=V_{1}^{2}+\partial_{z_{1}}V_{1}$ are small functions of $f$.
Substituting it into (\ref{3.3}) we have $$\varphi_{1}=V_{3}f^{2},$$
where $V_{3}=V_{1}\frac{\partial_{z_{1}}g}{g}-V_{2}-nV_{1}^{2}$ is a small function of $f$. Thus
$$\|\quad T(r,f)=S(r,f),$$
this is absurd.

Hence, we get $P\equiv0$. We rewritten (\ref{1.4}) as follows
$$f^{n}\partial_{z_{1}}f=g.$$
Set
$$\frac{g}{f^{n+1}}=a,$$
then
$$a=\frac{f^{n}\partial_{z_{1}}f}{f^{n+1}}=\frac{\partial_{z_{1}}f}{f}.$$
By using the lemma of logarithmic derivative, we know $\|\quad m(r,a)=S(r,f)$. Note that
$$\|\quad N(r,a)\leq N(r,\partial_{z_{1}}f)+N\left(r,\frac{1}{f}\right)\leq 2N(r,f)+\frac{1}{n}N\left(r,\frac{1}{g}\right)=S(r,f).$$
Thus,
$$\|\quad T(r,a)=S(r,f).$$

This completes the proof of Theorem \ref{thm1.4}.

\section{Proof of theorem 1.8}
\quad Set $P=P_{n-1}(f)$. If $g$ is a constant, then (\ref{1.20}) can be rewritten as
\begin{equation}\label{4.1}
f^{n}\partial_{z_{i}}f=P^{*}(f),
\end{equation}
where $P^{*}(f)$ is a differential polynomial of degree at most $n-1$ in $f$. Since (\ref{4.1}) satisfies the conditions of the generalised Clunie lemma, we obtain
$$\|\quad m(r,\partial_{z_{i}}f)=S(r,f), \quad m(r,f\partial_{z_{i}}f)=S(r,f).$$
If $\partial_{z_{i}}f\not\equiv0$, we get
$$\|\quad T(r,\partial_{z_{i}}f)=S(r,f), \quad T(r,f\partial_{z_{i}}f)=S(r,f),$$
in view of $f$ is an entire. Thus,
$$\|\quad T(r,f)\leq T(r,f\partial_{z_{i}}f)+T\left(r,\frac{1}{\partial_{z_{i}}f}\right)=S(r,f),$$
which yields a contradiction. If $\partial_{z_{i}}f\equiv0$, then $f=\alpha(z_{1},...,z_{i-1},z_{i+1},...,z_{m})$, which is of the from in Theorem \ref{thm1.20}. Thus, we assume that $g$ is not a constant in the following proof.

Suppose that $g$ is transcendental. Then the nonconstant entire function $f$ must be also transcendental, otherwise the left-hand side of (\ref{1.20}) would be a polynomial but the right-hand side of (\ref{1.20}) would be transcendental, which is
impossible. By Lemma \ref{lem2.6} yields that
$$\|\quad T(r,g(f))=T(r,f^{n}\partial_{z_{i}}f+P_{n-1}(f))=O(T(r,f)).$$
On the other hand, by Lemma \ref{lem2.5} yields that
$$\lim_{r\rightarrow\infty}\frac{T(r,g(f))}{T(r,f)}=\infty,$$
a contradiction. Thus, $g$ is a rational function.

We notice that $g$ has at most one pole in $\mathbb{C}$. Otherwise, by the Picard theorem, $f$ has at most one finite Picard value and thus would take one of the poles of $g$ infinitely many times, and then the right-hand side of (\ref{1.20})
would not be an entire function, a contradiction.

We distinguish two cases to discuss.

Case 1: If $g$ has no poles, then $g$ is a polynomial. We write (\ref{1.20}) as
\begin{equation}\label{4.2}
f^{n}\partial_{z_{i}}f+P=a_{k}f^{k}+a_{k-1}f^{k-1}+\cdots+a_{1}f+a_{0},
\end{equation}
where $a_{j}$, $j=0,1,...,k$ are distinct complex numbers, $k$ is an integer. It is easy to see that $k\leq n+1$, since by (\ref{4.2}),
\begin{align*}\label{4.3}
\|\quad km(r,f)&=m(r,f^{n}\partial_{z_{i}}f+P)+O(1)\nonumber\\
&\leq (n+1)m(r,f)+S(r,f).
\end{align*}

Next we discuss three subcases.

Subcase 1.1: $k\leq n-1$. We write (\ref{4.2}) as
\begin{equation*}\label{4.4}
f^{n}\partial_{z_{i}}f=a_{k}f^{k}+a_{k-1}f^{k-1}+\cdots+a_{1}f+a_{0}-P.
\end{equation*}
With $g$ is a constant similarly, we get a contradiction.

Subcase 1.2: $k=n$. We write (\ref{4.2}) as
\begin{equation}\label{4.6}
f^{n}(\partial_{z_{i}}f-a_{n})=a_{n-1}f^{n-1}+\cdots+a_{1}f+a_{0}-P.
\end{equation}
Since (\ref{4.6}) satisfies the conditions of the generalised Clunie lemma, we obtain
\begin{equation*}\label{4.7}
\|\quad m(r,\partial_{z_{i}}f-a_{n})=S(r,f), \quad m(r,f(\partial_{z_{i}}f-a_{n}))=S(r,f).
\end{equation*}
If $\partial_{z_{i}}f-a_{n}\not\equiv0$, we have
\begin{equation*}
\|\quad T(r,\partial_{z_{i}}f-a_{n})=S(r,f), \quad T(r,f(\partial_{z_{i}}f-a_{n}))=S(r,f),
\end{equation*}
in view of $f$ is an entire. Thus,
$$\|\quad T(r,f)\leq T(r,f(\partial_{z_{i}}f-a_{n}))+T\left(r,\frac{1}{\partial_{z_{i}}f-a_{n}}\right)=S(r,f),$$
which yields a contradiction. If $\partial_{z_{i}}f-a_{n}\equiv0$, we have $f=a_{n}z_{i}+\alpha(z_{1},...,z_{i-1},z_{i+1},...,z_{m})$, which is of the from in Theorem \ref{thm1.20}.

Subcase 1.3: $k=n+1$. We write (\ref{4.2}) as
\begin{equation*}\label{4.8}
f^{n}(\partial_{z_{i}}f-a_{n+1}f-a_{n})=a_{n-1}f^{n-1}+\cdots+a_{1}f+a_{0}-P.
\end{equation*}
If $\partial_{z_{i}}f-a_{n+1}f-a_{n}\not\equiv0$, we can get a contradiction by the same argument in Subcase 1.2. If $\partial_{z_{i}}f-a_{n+1}f-a_{n}\equiv0$, that is
$\partial_{z_{i}}f=a_{n+1}f+a_{n}$. By Lemma \ref{lem2.2}, we have $f=az_{i}+\alpha(z_{1},...,z_{i-1},z_{i+1},...,z_{m})$ or $f=a+e^{bz_{i}+\alpha(z_{1},...,z_{i-1},z_{i+1},...,z_{m})}$, which are of the from in Theorem \ref{thm1.20}.

Case 2: If $g$ has one pole $a\in\mathbb{C}$. We can then write $g$ in the following form:
$$g(w)=\frac{1}{(w-a)^{t}}p(w)$$
for $w\neq a$ in $\mathbb{C}$, where $t$ is a positive integer, and $p$ is a polynomial in $\mathbb{C}$ with $p(a)\neq0$. By (\ref{1.20}), $f(z)\neq a$ for any $z\in\mathbb{C}^{m}$. Hence we have that $f=a+e^{q(z)}$, $z\in\mathbb{C}^{m}$, where
$q$ is an entire function in $\mathbb{C}^{m}$, which is non-constant since $f$ is non-constant. By (\ref{1.20}) again, we obtain that
\begin{equation}\label{4.9}
(e^{q})^{n+t}\partial_{z_{i}}e^{q}+Q(e^{q})=p(a+e^{q(z)})=u(e^{q}),
\end{equation}
where $Q$ is a differential polynomial of degree at most $n+t$ in $e^{q}$, $u$ is an entire function in $\mathbb{C}$. By the same argument as above for proving $g$ to be rational, we have $u$ must be a polynomial. We can get the degree of $u$ at
most $n+t+1$ by the same argument in Case 1.

Next we discuss two subcases.

Subcase 2.1: $\deg u\leq n+t$. By the generalised Clunie lemma, we obtain
\begin{equation*}\label{4.10}
\|\quad m(r,\partial_{z_{i}}e^{q})=S(r,e^{q}).
\end{equation*}
If $\partial_{z_{i}}e^{q}\not\equiv0$, we have
\begin{equation*}
\|\quad T(r,\partial_{z_{i}}e^{q})=S(r,e^{q}),
\end{equation*}
in view of $q$ is an entire. Then $T(r,e^{q})=S(r,e^{q})$, a contradiction. If $\partial_{z_{i}}e^{q}\equiv0$, we have $\partial_{z_{i}}q\equiv0$, $f=a+e^{\alpha(z_{1},...,z_{i-1},z_{i+1},...,z_{m})}$, which is of the from in Theorem
\ref{thm1.20}.

Subcase 2.2: $\deg u=n+t+1$. We write (\ref{4.9}) as
\begin{equation*}\label{4.11}
e^{(n+t)q}(\partial_{z_{i}}e^{q}-c_{n+t+1}e^{q})=c_{n+t}e^{(n+t)q}+\cdots+c_{0}-Q(e^{q}),
\end{equation*}
where $c_{j}$, $j=0,1,...,n+t+1$ are distinct complex numbers. By the generalised Clunie lemma, we obtain
\begin{equation*}\label{4.12}
\|\quad m(r,\partial_{z_{i}}e^{q}-c_{n+t+1}e^{q})=S(r,e^{q}).
\end{equation*}
If $\partial_{z_{i}}e^{q}-c_{n+t+1}e^{q}\not\equiv0$, we have
\begin{equation*}
\|\quad T(r,\partial_{z_{i}}e^{q}-c_{n+t+1}e^{q})=T(r,(\partial_{z_{i}}q-c_{n+t+1})e^{q})=S(r,e^{q}),
\end{equation*}
in view of $q$ is an entire. Then $T(r,e^{q})=S(r,e^{q})$, a contradiction. If $\partial_{z_{i}}e^{q}-c_{n+t+1}e^{q}\equiv0$, we have $\partial_{z_{i}}q\equiv c_{n+t+1}$, $f=a+e^{c_{n+t+1}z_{i}+\alpha(z_{1},...,z_{i-1},z_{i+1},...,z_{m})}$, which
is of the from in Theorem \ref{thm1.20}.

This completes the proof of Theorem \ref{thm1.20}.

\section{Data Availability}
The data used to support the findings of this study are
included within the article.

\section{Conflict of interests}
The authors declare that there is no conflict of interests
regarding the publication of this paper.

\end{document}